%



\input amstex.tex

\magnification=\magstep1
\hsize=5.5truein
\vsize=9truein
\hoffset=0.5truein
\parindent=10pt
\newdimen\nagykoz
\newdimen\kiskoz
\nagykoz=7pt
\kiskoz=2pt
\parskip=\nagykoz
\baselineskip=12.7pt


\loadeufm \loadmsam \loadmsbm

\font\vastag=cmssbx10
\font\drot=cmssdc10
\font\vekony=cmss10
\font\vekonydolt=cmssi10
\font\cimbetu=cmssbx10 scaled \magstep1
\font\szerzobetu=cmss10

\font\scVIII=cmcsc8
\font\rmVIII=cmr8
\font\itVIII=cmti8
\font\bfVIII=cmbx8
\font\ttVIII=cmtt8

\def\cim#1{{\centerline{\cimbetu#1}}}
\def\szerzo#1{{\vskip0.3truein\centerline{\szerzobetu#1}}}
\def\alcim#1{{\medskip\centerline{\vastag#1}}}
\def\tetel#1#2{{{\drot#1}{\it\szukebb~#2\tagabb}}}
\long\def\biz#1#2{{{\vekony#1} #2}}
\def\kiemel#1{{\vekonydolt#1\/}}
\long\def\absztrakt#1#2{{\vskip0.4truein{\vekony#1} #2\vskip0.5truein}}
\def\szukebb{\parskip=\kiskoz}
\def\tagabb{\parskip=\nagykoz}
\def\vonal{{\vrule height 0.2pt depth 0.2pt width 0.5truein}}

\def\CC{{\Bbb C}}

\def\bmfd{{Banach manifold}}
\def\bsmfd{{Banach submanifold}}

\def\hsmfd{{Hilbert submanifold}}

\def\bvbdl{{Banach vector bundle}}

\def\Ker{\hbox{\rm Ker}}

\def\im{\hbox{\rm Im}}

\def\Oa{{\Omega}}
\def\oa{{\omega}}
\def\UU{{\frak U}}

\def\Hom{\text{\rm Hom}}
\def\End{\text{\rm End}}

\def\ad{{\text{\rm ad}}}

\def\vfld{{vector field}}
\def\recpr{{reciprocal pair}}
\def\genrecpr{{generalized \recpr}}
\def\lderv{{Lie derivative}}

\def\hspc{{Hilbert space}}

\def\vspc{{vector space}}
\def\tspc{{topological space}}

\def\cts{{continuous}}
\def\bdd{{bounded}}

\def\idml{{infinite dimensional}}

\def\idms{{infinite dimensions}}

\def\fdml{{finite dimensional}}

\def\fdms{{finite dimensions}}

\def\pscx{{pseudoconvex}}

\def\diam{{\hbox{\rm diam}}}
\def\homo{{homomorphism}}
\def\sbs{{Schauder basis}}
\def\ubs{{unconditional basis}}

\def\psh{{plurisubharmonic}}
\def\pshdom{{\psh\ domination}}

\def\st{{such that}}
\def\wrt{{with respect to}}

\def\delbar{{\bold{\bar\partial}}}

\def\<{{\langle}}
\def\>{{\rangle}}
\def\NN {{\Bbb N}}
\def\RR {{\Bbb R}}

\def\OO {{\Cal O}}

\def\LL {{\Cal L}}

\def\epsz{\varepsilon}
\def\fii{\varphi}
\def\fn{func\-tion}
\def\fns{func\-tions}
\def\holo{hol\-o\-mor\-phic}

\def\ulrg{{uniform local rate of growth}}

\def\mfd{manifold}
\def\smfd{submanifold}
\def\cpx{complex}
\def\cpt{compact}
\def\ses{short exact sequence}

\def\nbd{neighbor\-hood}

\def\bspc{Banach space}

\def\ka{{\kappa}}

\def\La{\Lambda}

\def\da{{\delta}}

\def\ga{{\gamma}}
\def\aa{{\alpha}}

\def\da{{\delta}}

\def\Ga{{\Gamma}}

\def\cln{{:}\;}

\def\Prop{Proposition}
\def\p#1{{\Prop~#1}}

\def\Th{{Theorem}}
\def\th{{theorem}}

\def\t#1{{\Th~#1}}

\newcount\minute    
\newcount\hour      
\newcount\hourMins  
%
%
\def\now%
{
%
  \minute=\time    
  \hour=\time \divide \hour by 60 
  \hourMins=\hour \multiply\hourMins by 60
  \advance\minute by -\hourMins 
  \zeroPadTwo{\the\hour}:\zeroPadTwo{\the\minute}%
}
%
%
\def\timestamp%
{
  \today\ \now
}
%
%
\def\today%
{
  \the\year-\zeroPadTwo{\the\month}-\zeroPadTwo{\the\day}%
}
%
%
\def\zeroPadTwo#1%
{
%
  \ifnum #1<10 0\fi    
  #1
}

{\ttVIII\jobname, \timestamp}
{\phantom.}
\vskip0.5truein
\cim{ON THE ANNIHILATOR OF A DOLBEAULT GROUP}
\szerzo{Imre Patyi\plainfootnote{${}^*$}{\rmVIII 
Supported in part by NSF grant 
DMS 0600059.}}
\absztrakt{ABSTRACT.}{We show that any Dolbeault cohomology group $H^{p,q}(D)$, 
	$p\ge0$, $q\ge1$,
	of an open subset $D$ of a closed finite codimensional complex
	Hilbert submanifold of $\ell_2$ is either zero or infinite
	dimensional.
	 We also show that any continuous character of the algebra
	of holomorphic functions of a closed complex Hilbert submanifold $M$
	of $\ell_2$ is induced by evaluation at a point of $M$.
	 Lastly, we prove that any closed split infinite dimensional complex
	Banach submanifold of $\ell_1$ admits a nowhere critical
	holomorphic function.




}

\def\sA{{1}}
\def\sB{{2}}
\def\sC{{3}}
\def\sD{{4}}
\def\sE{{5}}
\def\sF{{6}}
\def\sG{{7}}


\def\tBA{{\sB.1}}
\def\tBB{{\sB.2}}

\def\tCA{{\sC.1}}
\def\tCX{{\sC.2}}
\def\tCB{{\sC.3}}
\def\tCC{{\sC.4}}
\def\tCD{{\sC.5}}
\def\tCE{{\sC.6}}
\def\tCF{{\sC.7}}
\def\tCG{{\sC.8}}

\def\tDA{{\sD.1}}
\def\tDB{{\sD.2}}
\def\tDC{{\sD.3}}
\def\tDD{{\sD.4}}
\def\tDE{{\sD.5}}
\def\tDF{{\sD.6}}

\def\eDA{{\sD.1}}

\def\tEA{{\sE.1}}
\def\tEB{{\sE.2}}
\def\tEC{{\sE.3}}
\def\tED{{\sE.4}}

\def\tFA{{\sF.1}}
\def\tFB{{\sF.2}}

\def\tGA{{\sG.1}}
\def\tGB{{\sG.2}}
\def\tGC{{\sG.3}}
\def\tGD{{\sG.4}}

\def\rD{D}
\def\rDPV{DPV}
\def\rF{F}
\def\rLA{L1}
\def\rLB{L2}
\def\rLP{LP}
\def\rLf{Lf}
\def\rLT{LT}
\def\rM{M}
\def\rMc{Mc}
\def\rP{P}
\def\rPt{Pt}
\def\rS{S}
\def\rZ{Z}

\alcim{\sA. INTRODUCTION.}

	 In this paper we study ideals of \holo\ \fns\ and their connections with
	\holo\ \vfld{}s and Lie derivatives.
	 We generalize in \t\tFB\ a \th\ of Laufer on the \idml{}ity of
	certain Dolbeault groups.
	 We also generalize in \t\tDE\ a \th\ of Schottenloher on \cts\ characters
	of the algebra of \holo\ \fns\ of \cpx\ \bmfd{}s.
	 We show in \t\tDD\ that ideals of \holo\ \fns\ over a \cpx\ \bmfd{} without
	common zeros are often sequentially dense.
	 We prove that certain \cpx\ \bmfd{}s admit nowhere critical \holo\ \fns;
	see \Th{}s \tCC, \tCD, \tGD, \Prop{}s \tCE\ and \tCF.

\alcim{\sB. BACKGROUND.}

         In this section we collect some definitions and theorems that are useful for
        this paper.  
	 Some good sources of information on \cpx\ analysis on \bspc{}s are 
	[\rD, \rM, \rLA].

	 Put $B_X$ for the open unit ball of a \bspc\ $(X,\|\cdot\|)$,
	$\End(X)$ for the \bspc\ of \bdd\ linear operators $T{:}\;X\to X$ endowed
	with the operator norm, and $X^*$ for the dual space of $X$.
         A \kiemel{\cpx\ \bmfd} $M$ modelled on a \cpx\ \bspc\ $X$ is a paracompact
        Hausdorff space $M$ with an atlas of bi\holo{}ally related charts onto
        open subsets of $X$.
         A subset $N\subset M$ is called a \kiemel{closed \cpx\ \bsmfd} of $M$
        if $N$ is a closed subset of $M$ and for each point $x_0\in N$ there are
        an open \nbd\ $U$ of $x_0$ in $M$ and a bi\holo\ map $\fii{:}\;U\to B_X$ onto
        the unit ball $B_X$ of $X$ that maps the pair $(U,U\cap N)$ to
        a pair $(B_X,B_X\cap Y)$ for a closed \cpx\ linear subspace $Y$ of $X$.
         The \smfd\ $N$ is called a \kiemel{split} or \kiemel{direct} \bsmfd\ of $M$
        if at each point $x_0\in N$ the corresponding subspace $Y$ has a direct complement
        in $X$.
	 Following [\rLB] by Lempert we say that \kiemel{\pshdom} is possible in
	a complex \bmfd\ $M$ if for any $u{:}\;M\to\RR$ locally upper bounded there is
	a $\psi{:}\;M\to\RR$ \cts\ and \psh\ \st\ $u(x)<\psi(x)$ for all $x\in M$.
	 This is a kind of \holo\ convexity property of $M$.

\tetel{\t\tBA.}{{\rm(Lempert, [\rLB])}
	 If $X$ is a \bspc\ with an \ubs\ and\/ $\Oa\subset X$ is \pscx\ open,
	then \pshdom\ is possible in\/ $\Oa$.
}

	 We make use of the following vanishing \t\tBB.

\tetel{\t\tBB.}{Let $X$ be a \bspc\ with a \sbs, $\Oa\subset X$ \pscx\ open,
	$M\subset\Oa$ a closed split \cpx\ \bsmfd\ of\/ $\Oa$, and $E\to M$ a
	\holo\ \bvbdl.
	 Suppose that \pshdom\ is possible in every \pscx\ open subset of\/ $\Oa$.
	 Then the following hold.
\vskip0pt
	{\rm(a)} Let $\OO^E\to M$ be the sheaf of germs of \holo\ sections of
	$E\to M$.
	 Then the sheaf cohomology group $H^q(M,\OO^E)$ vanishes for all $q\ge1$.
\vskip0pt
	{\rm(b)} Any \holo\ \fn\ $f\in\OO(M,Z)$ into any \bspc\ $Z$ can be extended
	to an $\tilde f\in\OO(\Oa,Z)$ with $\tilde f(x)=f(x)$ for all $x\in M$.
\vskip0pt
	{\rm(c)} If\/ $0\to E'\to E\to E''\to0$ is a pointwise split \ses\ of \holo\
	\bvbdl{}s over $M$, then it admits a \holo\ global splitting over $M$.
\vskip0pt
	{\rm(d)} If $(f_n)\in\OO(M,\ell_2)$ is nowhere zero on $M$, then there is
	a $(g_n)\in\OO(M,\ell_2)$ with\/ $\sum_{n=1}^\infty f_n(x)g_n(x)=1$
	for all $x\in M$, where the series converges absolutely and uniformly
	on every \cpt\ subset of $M$.
}

\biz{Proof.}{These are special cases of the vanishing \th\ of [\rLP]; for (d)
	see also [\rDPV].
}

\alcim{\sC. RECIPROCAL PAIRS.}

	 In this section we look at the following simple notion.
	 Let $M$ be a \cpx\ \bmfd, $f\in\OO(M)$ a \holo\ \fn, and
	$v\in\OO(M,T^{1,0}M)$ a \holo\ (tangent) \vfld\ on $M$.
	 We call $f,v$ a \kiemel{\recpr} if the \lderv\
	$(\LL_vf)(x)=(vf)(x)=(df)(x)v(x)=1$ for all $x\in M$.
	 If $f,v$ is a \recpr\ on $M$, then the Fr\'echet differential
	$df$ and $v$ do not have zeros on $M$; in particular,
	$f$ is nowhere critical on $M$, i.e., it has no points $x\in M$
	with $(df)(x)=0$.
	 A useful weakening of the notion of a \recpr\ is a \genrecpr\ defined as follows.
	 Let $n\ge1$, $f_i\in\OO(M)$ and $v_i\in\OO(M,T^{1,0}M)$ for $i=1,\ldots,n$.
	 We call $f_1,\ldots,f_n$, $v_1,\ldots,v_n$ a \kiemel{\genrecpr} if
	$\sum_{i=1}^nv_if_i=1$ on $M$.
	 The simple proof of \p\tCA\ below is omitted.

\tetel{\p\tCA.}{Let $M$ be a \cpx\ \bmfd.
\vskip0pt
	{\rm(a)} If $M$ admits a \genrecpr, and $M'$ is bi\holo\ to $M$, then
	$M'$ also admits an analogous \genrecpr.
\vskip0pt
	{\rm(b)} If $D\subset M$ is open, and $M$ admits a \genrecpr\ $f_i,v_i$ for
	$i=1,\ldots,n$, then
	$D$ admits the \genrecpr\ $f_i|D,v_i|D$ of restrictions.
\vskip0pt
	{\rm(c)} Suppose that $M$ admits \fns\ $f_i\in\OO(M)$ for $i=1,\ldots,n$
	without common critical points, and look at the \ses\
$$
	0\to K\to (T^{1,0}M)^n\to M\times\CC\to0
 $$
 	of \holo\ \bvbdl{}s over $M$, where the third mapping is\/ 
	$(\xi_i)\mapsto\sum_{i=1}^n(df_i)\xi_i$,
	and the second mapping is inclusion; so $K$ is the kernel of the third mapping.
	 Then there are \holo\ \vfld{}s $v_i\in\OO(M,T^{1,0}M)$ for $i=1,\ldots,n$
	with $\sum_{i=1}^nv_if_i=1$ on $M$
	if and only if the above \ses\ of \holo\ \bvbdl{}s splits \holo{}ally over $M$.
	 The latter is the case if the sheaf cohomology group $H^1(M,\OO^K)$ vanishes,
	where $\OO^K\to M$ is the sheaf of germs of \holo\ sections of the \holo\ \bvbdl\
	$K\to M$.
\vskip0pt
	{\rm(d)} If $M'$ is a \cpx\ \bmfd\ that admits a \genrecpr\ $f'_i,v'_i$ 
	for $i=1,\ldots,n$, and
	$M''$ is any \cpx\ \bmfd, then $M=M'\times M''$ also admits an analogous
	\recpr\ $f_i,v_i$ given by $f_i(x',x'')=f'_i(x')$ and $v_i(x',x'')=(v'_i(x'),0)$
	for $i=1,\ldots,n$.
}

\tetel{\t\tCX.}{Let $X$ be a \bspc\ with a \sbs, $\Oa\subset X$ \pscx\ open,
	$M\subset\Oa$ a closed split \cpx\ \bsmfd\ of\/ $\Oa$, and suppose that
	\pshdom\ is possible in every \pscx\ open subset of\/ $\Oa$ (the last is
	guaranteed by Lempert's \t\tBA\ if $X$ has an \ubs).
	 If $f_i\in\OO(M)$, $i=1,\ldots,n$, have no common critical points in $M$, 
	then there are $v_i\in\OO(M,T^{1,0}M)$, $i=1,\ldots,n$,
	\st\/ $\sum_{i=1}^nv_if_i=1$ on $M$.
}

\biz{Proof.}{\t\tBB(a) implies that \p\tCA(c) applies,
	completing the proof of \t\tCX.
}

	We recall a deep theorem of Forstneri\v{c}.

\tetel{\t\tCB.}{{\rm(a) (Forstneri\v{c}, [\rF])}
	Every Stein manifold admits a nowhere critical \holo\ \fn.
\vskip0pt
	{\rm(b)} Every Stein manifold without isolated points admits a \recpr.
}

\biz{Proof of (b).}{As Theorem~B of Cartan, Oka, Serre implies the vanishing of
	the relevant cohomology group, (b) follows from (a) via \p\tCA(c).
}

\tetel{\t\tCC.}{If $M$ is a \cpx\ \bmfd\ satisfying\/ {\rm(a)} or\/ {\rm(b)} below,
	and $D\subset M$ is any open subset, then $D$ admits a \recpr.
\vskip0pt
	{\rm(a)} $M$ is any \bspc.
\vskip0pt
	{\rm(b)} $M$ is the $W^{(k)}_2$-Sobolev space of mappings $x{:}\;K\to N$ that have
	$k$ derivatives in $L_2$, where $K$ is any \cpt\ smooth \mfd, $N$ is any
	Stein \mfd, and $k$ is any integer with $2k>\dim_{\RR}(K)$.
}

\biz{Proof.}{(a) Choose a linear functional $f\in M^*$ and a vector $v\in M$ \st\
	$f(v)=1$.
	 Regard $v$ as a constant \vfld\ on $M$.
	 (b) Forstneri\v{c}'s \t\tCB(a) gives a nowhere critical $g\in\OO(N)$,
	choose any probability Radon measure $\mu$ on $K$, e.g., $\mu=\da_{t_0}$
	the Dirac delta measure concentrated at a point $t_0\in K$,
	and define $f\in\OO(M)$ by $f(x)=\int_{t\in K}g(x(t))\,d\mu(t)$ for $x\in M$.
	 It is easy to check that $f$ is indeed a \holo\ \fn\ on $M$, and
	that it has no critical points in $M$.
	 \t\tBB(a) implies that \p\tCA(c) applies and
	gives us a \vfld\ $v\in\OO(M,T^{1,0}M)$ with $vf=1$.
	 \p\tCA(b) concludes the proof of \t\tCC.
}

	 We omit the simple proof of the following \t\tCD, cf.\ \t\tGD.

\tetel{\t\tCD.}{If $X$ is a separable \bspc, its dual $X^*$ is nonseparable,
	$\Oa\subset X$ open, and $M\subset\Oa$ a closed \cpx\
	\bsmfd\ of\/ $\Oa$ of finite codimension,
	then there is a linear functional $\xi\in X^*$ whose restriction $f=\xi|M$
	is a nowhere critical \fn\ $f\in\OO(M)$ on $M$.
	 Further, if $X$ has an \ubs\ and\/ $\Oa$ is \pscx, 
	then there is a $v\in\OO(M,T^{1,0}M)$ with $vf=1$ on $M$.
}

	 In \t\tCD\ the \bspc\ $X$ can be $X=\ell_1\times Y$, where
	$Y$ is any \bspc\ with an \ubs.
	 Sometimes a \recpr\ can be constructed with explicit computation as in
	\p\tCE\ below, whose easy proof we omit.

\tetel{\p\tCE.}{Let $X$ be a \fdml\ or separable \hspc\ of dimension at least two
	with standard coordinate \fns\ $x_1,x_2,\ldots$, and $M\subset X$ the smooth 
	hypersurface defined by $g(x)=0$, where $g(x)=-1+\sum_jx_j^2$.
	 Then a \recpr\ $f,v$ on $M$ is given by
	$f(x)=x_1+ix_2$, and
	$v(x)=-ix_1^2D_2+ix_1x_2D_1+D_1-x_1E$,
	where $i=\sqrt{-1}$, $D_j=\frac{\partial}{\partial x_j}$ is the usual Wirtinger
	derivative \wrt\ $x_j$, and $E$ is the Euler derivative $(Eh)(x)=(dh)(x)x$
	for $h\in\OO(X)$, i.e., $E=\sum_j x_jD_j$.
}

	 A generalization of part of \p\tCE\ is given below in \p\tCF,
	whose proof is clear from the definitions.

\tetel{\p\tCF.}{Let $X=X'\times X''$ be a direct decomposition of \bspc{}s,
	$g_1\in\OO(X')$, $g_2\in\OO(X'')$ entire \fns, and suppose that
	the hypersurface $M'\subset X'$ defined by\/ $0=g_1(x')$
	is smooth, $(dg_1)(x')\not=0$ for $x'\in M'$, there is an $f_1\in\OO(X')$ \st\
	$f_1$ is nowhere critical on $X'$ and so is $f_1|M'$ on $M'$, $g_2(0)=0$,
	and the only critical point of $g_2$ in $X''$, if any, is $x''=0$.
	 Then the hypersurface $M\subset X$ defined by\/ $0=g(x',x'')$,
	where $g(x',x'')=g_1(x')+g_2(x'')$, is smooth and the \fn\ $f\in\OO(X)$
	defined by $f(x',x'')=f_1(x')$ is nowhere critical on $M$.
}

	 \p\tCE\ (without the vector field $v$) is a special case of \p\tCF,
	where $X'=\CC^2$, $g_1(x'_1,x'_2)=-1+(x'_1)^2+(x'_2)^2$, 
	$f_1(x'_1,x'_2)=x'_1+ix'_2$, and $g_2(x''_1,x''_2,\ldots)=\sum_j(x''_j)^2$.

	 We can also make various other special cases of \p\tCF, e.g., $X',g_1,f_1$
	as above, $X''=\ell_p$, $1\le p\le\infty$, and
	$g_2(x'')=\sum_{j=1}^\infty a_j(x''_j)^{r_j}$, where the $r_j\ge1$ are
	integers, and the coefficients $0\not=a_j\in\CC$ are such that $g_2\in\OO(X'')$,
	e.g., $a_j\to0$ fast enough as $j\to\infty$.

	 \p\tCG\ below is obvious from elementary linear algebra.

\tetel{\p\tCG.}{Let $M'$ be a \cpx\ \bmfd, $M\subset M'$ a closed \cpx\ \bsmfd\ of $M'$
	of finite codimension $k\ge1$.
	 If $f_\ka\in\OO(M')$ for $\ka=0,\ldots,k$ satisfy that $(df_\ka)(x)$ for 
	$\ka=0,\ldots,k$ are linearly independent at each point $x\in M$, then the
	restrictions $f_\ka|M\in\OO(M)$ for $\ka=0,\ldots,k$ have no common critical
	points in $M$.
}

	 If $M'$ is a \bspc\ of dimension at least $k+1$, then we can choose the
	$f_\ka\in(M')^*$ to be any $k+1$ linearly independent linear functionals in \p\tCG.

\alcim{\sD. FAMILIES OF FUNCTIONS.}

	 In this section we look at the following notion.
	 Let $M$ be a topological space, $f_n{:}\;M\to\CC$, $n\ge1$, a sequence of
	numerical \fns.
	 We say that the sequence $(f_n)$ has a \kiemel{\ulrg} (as $n\to\infty$)
	if there are a sequence of constants $L_n>0$, $n\ge1$, an open covering $\UU$
	of $M$, and a \fn\ $C{:}\;\UU\to(0,\infty)$ \st\ if $U\in\UU$, $x\in U$, and $n\ge1$,
	then $|f_n(x)|\le C(U)L_n$.

	 In other words, near each point $x_0\in M$ our \fns\ $f_n(x)$ are \bdd,
	and they grow no faster than $(\text{const})L_n$ for a sequence of
	constants $L_n$, $n\ge1$.

	 We omit the proof of the following obvious \p\tDA.

\tetel{\p\tDA.}{{\rm(a)} If $M$ is a Lindel\"of space, then
	any sequence of locally bounded \fns\ $f_n{:}\;M\to\CC$, $n\ge1$, has a \ulrg.
\vskip0pt
	{\rm(b)} Let $M$ be a para\cpt\ Hausdorff space, and $f_n{:}\;M\to\CC$, $n\ge1$,
	a sequence of \fns.
	 Then\/ $(f_n)$ has a \ulrg\/ $(L_n)$ if and only if there is a \cts\ \fn\
	$\ga{:}\;M\to(0,\infty)$ with\/ $|f_n(x)|\le\ga(x) L_n$ for all $x\in M$ and
	$n\ge1$.
}

	 On a \fdml\ \cpx\ \mfd\ any sequence of \holo\ \fns\ has a \ulrg\ by \p\tDA(a),
	unlike on an \idml\ one.
	 We now show that certain natural sequences do.

\tetel{\p\tDB.}{{\rm(a)} Let $X$ be a \cpx\ \bspc, $B_X$ its open unit ball, and
	$f\in\OO(B_X)$.
	 If there is a bound\/ $0\le M<\infty$ \st\/ $|f(x)|\le M$ for\/ $\|x\|<1$, then\/
	$|f^{(n)}(x)\xi_1\ldots\xi_n|\le\frac{Mn^n}{(1-\|x\|)^n}\|\xi_1\|\ldots\|\xi_n\|$
	for $x\in B_X$, $\xi_1,\ldots,\xi_n\in X$, and $n\ge0$.
\vskip0pt
	{\rm(b)} Let $X$ be a \cpx\ \bspc, $\Oa\subset X$ open, $f\in\OO(\Oa)$,
	and $\xi^{(n)}_j\in\overline{B_X}$ for\/ $1\le j\le n$, $n\ge0$.
	 Then the sequence of \fns\ $f^{(n)}(x)\xi^{(n)}_1\ldots\xi^{(n)}_n$
	for $x\in\Oa$, $n\ge0$, has a \ulrg\ $n^{2n}$.
}

\biz{Proof.}{Part (a) follows from the usual Cauchy estimate for a polydisc.
	Part (b) follows on applying (a) locally on balls in $\Oa$ on which $f$ is \bdd.
}

\tetel{\p\tDC.}{If $X$ is a \bspc, $\Oa\subset X$ open, and $f\in\OO(\Oa)$, then the
	following hold.
\vskip0pt
	{\rm(a)} If $v_n\in\OO(\Oa,X)$, $n\ge1$, is a sequence of \holo\ \vfld{}s 
	that has a \ulrg\ of\/ $1$, then the sequence of Lie derivatives
	 $f_n=\LL_{v_n}\LL_{v_{n-1}}\ldots\LL_{v_1}f\in\OO(\Oa)$ has a \ulrg\ of
	$n^{5n^2}$.
\vskip0pt
	{\rm(b)} Let $P\in\OO(\Oa,\End(X))$ be a \holo\ \fn\ with operator values,
	$\xi_n\in\overline{B_X}$, $n\ge1$, and $N$ the set of all finite sequences
	$n=(n_1,\ldots,n_s)$ of natural numbers $n_j\ge1$ for $s\ge1$.
	 Define $f_n\in\OO(\Oa)$, $n\in N$, by
	 $f_{n_1\ldots n_s}=\LL_{P\xi_{n_1}}\ldots\LL_{P\xi_{n_s}}f$.
	 Then\/ $(f_n)$ has a \ulrg\ of $s^{5s^2}$, where $s$ is the length of $n\in N$,
	i.e., there are an open covering\/ $\UU$ of\/ $\Oa$, and a \fn\ $C{:}\;\UU\to(0,\infty)$
	\st\ if $U\in\UU$, $x\in U$, and $n\in N$, 
	then\/ $|f_{n_1\ldots n_s}(x)|\le C(U)s^{5s^2}$.
}

\biz{Proof.}{(a) The \fn\ $f_n$ is the sum of $n!$ products, whose $n+1$ factors of each are
	at most $n$th derivatives of $f$ and at most $(n-1)$st derivatives of the $v_j$.
	 Each factor has a \ulrg\ of $n^{2n}$ by \p\tDB(b), each product of $n+1$ of
	them $(n^{2n})^{n+1}$, and the sum $f_n$ of the $n!$ such products
	$n!(n^{2n})^{n+1}\le n^{5n^2}$.
	 (The accumulation of the constants is much less severe than the growth of
	$n^{5n^2}$, hence it can be incorporated in the said growth rate.)

	 As (b) follows just in the same way as (a) does, we omit the rest of
	the proof of \p\tDC.
}

\tetel{\t\tDD.}{Let $X$ be a \bspc\ with a \sbs, $\Oa\subset X$ \pscx\ open, $M\subset\Oa$
	a closed split \cpx\ \bsmfd\ of\/ $\Oa$, $N$ a countable set, and
	$f_n\in\OO(M)$, $n\in N$.
	 Suppose that \pshdom\ is possible is every \pscx\ open subset of $\Oa$.
	 If the \fns\ $f_n$, $n\in N$, have no common zeros, and have a \ulrg,
	then there are \holo\ \fns\ $g_n\in\OO(M)$, $n\in N$, \st\
	$\sum_{n\in N} f_n(x)g_n(x)=1$ for all $x\in M$, where the series converges
	absolutely and uniformly on every \cpt\ subset of $M$.
}

\biz{Proof.}{There are an open covering $\UU$ of $M$, a \fn\ $C{:}\;\UU\to(0,\infty)$,
	and constants $L_n>1$, $n\in N$, \st\ if $U\in\UU$, $x\in U$, and $n\in N$, then
	$|f_n(x)|\le C(U)L_n$.
	 Let $i{:}\;N\to\{1,2,3,\ldots\}$ be an injection, $L'_n=2^{i(n)}L_n$, $n\in N$,
	$H=\{z=(z_n)_{n\in N}{:}\;z_n\in\CC,\|z\|=(\sum_{n\in N}|z_n|^2)^{1/2}<\infty\}$ our \hspc,
	and $F\in\OO(M,H)$ defined by $F(x)=(F_n(x))$, where $F_n(x)=f_n(x)/L'_n$, $n\in N$.
	 Then $F$ is indeed a \holo\ \fn\ $F{:}\;M\to H$, and $F(x)\not=0$ for every $x\in M$.
	 \t\tBB(d) applies and gives a \holo\ $G=(G_n)\in\OO(M,H)$
	with $1=F(x)\cdot G(x)=\sum_{n\in N}F_n(x)G_n(x)$ for $x\in M$, where the convergence
	is absolute and uniform on every \cpt\ subset of $M$.
	 Letting $g_n=G_n/L'_n$, $n\in N$, completes the proof of \t\tDD.
}

	 In the setting of \t\tDD\ suppose that $I\subset\OO(M)$ is an ideal 
	of $\OO(M)$, and $I$ is 
	sequentially closed in the sense that if $f(x)=\sum_{n=1}^\infty f_n(x)g_n(x)$ 
	for $x\in M$, where $f_n\in I$, $g_n\in\OO(M)$, $n\ge1$, and the series
	converges absolutely and uniformly on every \cpt\ subset of $M$,
	then $f\in\OO(M)$ also lies in the ideal $I$.
	 Then $I=(1)$ is the unit ideal if and only if $I$ admits a sequence $f_n\in I$,
	$n\ge1$, without common zeros and with a \ulrg.

\tetel{\t\tDE.}{{\rm(a)} Let $X$ be a \bspc\ with a \sbs, $M\subset X$ a closed split
	\cpx\ \bsmfd\ of $X$, and suppose that \pshdom\ is possible in every \pscx\ open
	subset of $X$.
	 If $\chi{:}\;\OO(M)\to\CC$ is a \cts\ character of the algebra $\OO(M)$ (i.e.,
	$\chi$ is a multiplicative linear functional, $\chi(1)=1$, and
	there is a \cpt\ set $K\subset X$ with\/ $|\chi(f)|\le\sup_{x\in K}|f(x)|$
	for all $f\in\OO(M)$), then there is a point $x_0\in M$ with $\chi(f)=f(x_0)$
	for all $f\in\OO(M)$.
\vskip0pt
	{\rm(b)} Let $X$ be a \bspc\ with an \ubs, $\Oa\subset X$ \pscx\ open,
	$M\subset\Oa$ a closed split \cpx\ \bsmfd\ of\/ $\Oa$, and $\chi{:}\;\OO(M)\to\CC$
	a \cts\ character, then there is a point $x_0\in M$ with $\chi(f)=f(x_0)$
	for all $f\in\OO(M)$.
}

\biz{Proof.}{(a) Let $\xi_n\in X^*$, $n\ge1$, be the coordinate functionals of a bimonotone
	\sbs\ of $X$.
	 Thus there is a bound $0\le B<\infty$ with
	$|\xi_n(x)|\le B\|x\|$ for $x\in X$ and $n\ge1$. 
	 Look at the \fns\ $f_n\in\OO(X)$ defined by $f_n=\xi_n-\chi(\xi_n|M)$ for $n\ge1$.
	 Note that $|f_n(x)|\le B\|x\|+B\diam(K\cup\{0\})$ for $x\in X$ and $n\ge1$.
	 Hence the sequence $f_n$, $n\ge1$, has a \ulrg.
	 Let $I=\Ker\,\chi$, and note that $f_n|M\in I$ for $n\ge1$.

	 If $f_n|M$, $n\ge1$, have no common zeros in $M$, then by \t\tDD\ and 
	the remark following its proof we find that $I=(1)$, i.e., $\chi=0$, 
	which contradicts that $\chi(1)=1$.

	 Hence there is a point $x_0\in M$ with $f_n(x_0)=0$, i.e., 
$$
	f(x_0)=\chi(f|M)
\tag\eDA
 $$
	for $f=\xi_n$ for all $n\ge1$.
	 As $(\eDA)$ subsists for any polynomial of finitely many $\xi_n$,
	and as any $f\in\OO(X)$ is the limit of a sequence of such polynomials
	uniformly on $K\cup\{x_0\}$, we see that $(\eDA)$ holds for all $f\in\OO(X)$.
	 As any $f\in\OO(M)$ can be extended to a \holo\ $\tilde f\in\OO(X)$
	with $f=\tilde f|M$ by \t\tBB(b), the proof of (a)
	is complete.
	 Since (b) follows from (a) and from Zerhusen's embedding \th\ [\rZ]
	by embedding $M$ as a closed split \cpx\ \bsmfd\ $M'$ of a \bspc\ $X'$ with a \ubs,
	the proof of \t\tDE\ is complete.
}

	 It was shown much earlier by Schottenloher [\rS], see also [\rM],
	that a \cts\ character of $\OO(M)$ is a point evaluation of $M$
	if $M$ is a Riemann domain spread over a \pscx\ open subset
	of a \bspc\ with a \sbs.
	 It is unclear whether the statement of \t\tDE\ follows from the above mentioned
	result of Schottenloher.
	 It is, however, possible to replace a part of the proof of \t\tDE(a) 
	by an application of his result.
	 Indeed, one can consider the character $\chi'$ of $\OO(X)$ defined by 
	$\chi'(f)=\chi(f|M)$, apply his result to find a point $x_0\in X$ with $(\eDA)$
	for $f\in\OO(X)$, and conclude as above by invoking the extension \t\tBB(b)
	from $M$ to $X$ and from $M\cup\{x_0\}$ to $X$, should $x_0$ lie outside $M$.

	 Note that if $M'$ is a \holo\ covering \bmfd\ with countably many leaves of 
	an $M$ as in \t\tDE(b), then
	$M'$ is bi\holo\ to a closed split \cpx\ \bsmfd\ $M''$ of a \bspc\ with an \ubs,
	hence any \cts\ character of $\OO(M')$ is gotten by evaluation at a point of $M'$.

	 \t\tDF\ below follows from standard linear algebra and \t\tDE.

\tetel{\p\tDF.}{Let $M$ be a \cpx\ \bmfd\ as in \t\tDE, $E=\End(\CC^n)$ the
	algebra of \cpx\ $n$ by $n$ matrices, $n\ge1$, $A{:}\;\OO(M)\to E$ a \cpx\ algebra
	\homo\ with $A(1)=1$, and $I=\Ker\,A$.
	 Choose a basis of\/ $\CC^n$ so that the commuting matrices $A(f)$, $f\in\OO(M)$,
	are simultaneously upper triangular \wrt\ the chosen basis.
	 So $A(f)=[a_{ij}(f)]_{ij=1}^n$, and $a_{ij}(f)=0$ for $i>j$.
	 If at least one of the characters $a_{ii}{:}\;\OO(M)\to\CC$, $i=1,\ldots,n$, is \cts,
	then the ideal $I$ has a common zero.
}

\alcim{\sE. LIE DERIVATIVES AND IDEALS.}

	 In this section we look at ideals of \holo\ \fns\ that relate to Lie derivatives.

\tetel{\p\tEA.}{Let $X$ be a \bspc, $\Oa\subset X$ open, $f\in\OO(\Oa)$,
	$v_1,\ldots,v_n\in\OO(\Oa,X)$ \holo\ \vfld{}s, $n\ge1$, and $x_0\in\Oa$.
	 If $f$ has a zero at least of order $n$ at the point $x_0$ (i.e., $f^{(i)}(x_0)=0$
	for $i=0,\ldots,n-1$), then $(\LL_{v_1}\ldots\LL_{v_n}f)(x_0)=
	f^{(n)}(x_0)v_1(x_0)\ldots v_n(x_0)$ holds for the iterated Lie derivative.
}

\biz{Proof.}{This follows from the rules of differentiation such as the product rule.
}

\tetel{\p\tEB.}{Let $N$ be as in \p\tDC(b), $M$ a connected \cpx\ \bmfd, $f\in\OO(M)$, and 
	$v_n\in\OO(M,T^{1,0}M)$ \holo\ \vfld{}s for $n\ge1$.
	 If the set of values $v_n(x_0)\in T^{1,0}_{x_0}M$, $n\ge1$, is dense in a
	\nbd\ of\/ $0$ in the \bspc\ $T^{1,0}_{x_0}M$ at a point $x_0\in M$, and
	$(\LL_{v_{n_1}}\ldots\LL_{v_{n_s}}f)(x_0)=0$ for all $n\in N$, then $f$
	is a constant $f(x_0)$ on $M$.
}

\biz{Proof.}{Let $s$ be the vanishing order of $f-f(x_0)$ at $x_0$.
	 If $s=\infty$, then we are done.
	 If $1\le s<\infty$, then we conclude the proof of \p\tEB\ by an application of
	\p\tEA\ to the restrictions to an open \nbd\ of $x_0$ bi\holo\ to an open set $\Oa$
	in a \bspc\ $X$.
}

\tetel{\p\tEC.}{{\rm(a)} Let $X$ be a \bspc, and $P\in\End(X)$ a projection, i.e., $P^2=P$.
	 If $x_n\in X$, $n\ge1$, is dense in the unit ball $B_X$, then\/ $\{Px_n{:}\; n\ge1\}$
	has a subset contained and dense in the ball $\frac{1}{\|P\|}B_{PX}$ of 
	the image \bspc\ $PX$.
\vskip0pt
	{\rm(b)} Let $N$ be as in \p\tDC(b), $X$ a separable \bspc,
	$\xi_n\in X$, $n\ge1$, dense in $B_X$,
	$\Oa\subset X$ open, $M\subset\Oa$ a connected closed split \cpx\ \bsmfd\ of\/ $\Oa$,
	$P\in\OO(M,\End(X))$ a \holo\ operator \fn\ with $P(x)P(x)=P(x)$ and\/
	$\im\,P(x)=P(x)X=T_xM$ for all $x\in M$, $f\in\OO(M)$, and 
	$f_n=\LL_{P\xi_{n_1}}\ldots\LL_{P\xi_{n_s}}f\in\OO(M)$ for $n\in N$.
	 If the \fns\ $f_n$, $n\in N$, have a common zero $x_0\in M$, then
	$f$ is a constant $f(x_0)$.
}

\biz{Proof.}{(a) Fix $\epsz>0$, $x_0\in PX$ with $\|x_0\|<\frac{1}{\|P\|}$, and choose
	an $\eta>0$ so small that $\|x_0\|+\eta\|P\|<\frac{1}{\|P\|}$ and $\|P\|\eta<\epsz$.
	 There is an $n\ge1$ with $\|x_n-x_0\|<\eta$.
	 Then $Px_0=x_0$, $\|Px_n-Px_0\|<\|P\|\eta<\epsz$, and
	$\|Px_n\|\le\|x_0\|+\|P\|\eta<\frac{1}{\|P\|}$.

	(b) As the \vfld{}s $v_n(x)=P(x)\xi_n$, $n\ge1$,  have values dense near zero in
	$T_xM$ by (a) for each $x\in M$, an application of \p\tEB\ completes the proof
	of \p\tEC.
}

\tetel{\t\tED.}{Let $X$ be a \bspc\ with a \sbs, $\Oa\subset X$ \pscx\ open,
	$M\subset\Oa$ a connected closed split \cpx\ \bsmfd\ of\/ $\Oa$, and $f_0\in\OO(M)$.
	 Suppose that \pshdom\ is possible in every \pscx\ open subset of\/ $\Oa$.
	 If $f_0$ is not constant zero on $M$, then there are iterated Lie derivatives
	$f_n\in\OO(M)$, $n\ge1$, of $f_0$, and \holo\ \fns\ $g_n\in\OO(M)$, $n\ge0$, with
	$\sum_{n=0}^\infty f_n(x)g_n(x)=1$ for all $x\in M$, where the series converges
	absolutely and uniformly on every \cpt\ subset of $M$.
	 Further, if $I\subset\OO(M)$ is a nonzero ideal that is closed under 
	Lie derivation (i.e., $vf\in I$ for $f\in I$ and $v\in\OO(M,T^{1,0}M)$), then
	$I$ is sequentially dense in $\OO(M)$.
}

\biz{Proof.}{If $f_0$ is constant on $M$, then letting $g_0=1/f_0$ we are done.
	  Suppose now that $f_0$ is not constant on $M$, and extend $f_0$ from $M$
	 to $f_0\in\OO(\Oa)$ by \t\tBB(b).
	  As $M$ is a split \bsmfd\ of $\Oa$, locally the trivial bundle $M\times X$
	 splits as $TM\oplus E$.
	  By \t\tBB(c) there is a global splitting, i.e., we can
	 write $M\times X=TM\oplus E$, where $E\to M$ is a \holo\ Banach vector subbundle
	 of $M\times X$.
	  Define $P'\in\OO(M,\End(X))$ by projecting $(x,\xi)\in M\times X$ to
	 $P'(x)\xi\in T_xM$ in the above global direct decomposition of $M\times X$.
	  \t\tBB(b) gives us a \holo\ extension $P\in\OO(\Oa,\End(X))$
	 with $P|M=P'$.
	  Choose a sequence $\xi_n\in X$, $n\ge1$, dense in $B_X$, and define
	 $f_n\in\OO(\Oa)$ for $n\in N$, where $N$ is as in \p\tDC(b), by
	 $f_n=\LL_{P\xi_{n_1}}\ldots\LL_{P\xi_{n_s}}f_0$.
	  Let $N'=N\cup\{0\}$, and note that the \fns\ $f_n\in\OO(\Oa)$, $n\in N'$,
	 have a \ulrg\ by \p\tDC(b).
	  Hence $f_n|M\in\OO(M)$, $n\in N'$, also has a \ulrg\ on $M$, and no common zeros
	 in $M$ by \p\tEC(b).
	  \t\tDD\ applies and completes the proof of \t\tED.
}

	 An example of a proper dense ideal $I\subset\OO(\CC)$ that is closed under (Lie)
	derivation is $I=\{f\in\OO(\CC){:}\;\text{ord}_{n}(f)\to\infty
	\text{\ as\ }n\to\infty\text{\ in\ }\NN\}$,
	where $\text{ord}_{z_0}(f)$ is the vanishing order of $f$ at the point $z_0\in\CC$.

\alcim{\sF. THE ANNIHILATOR OF A DOLBEAULT GROUP.}

	 In this section we generalize to certain \idml\ \cpx\ \bmfd{}s the following
	\t\tFA\ of Laufer.

\tetel{\t\tFA.}{{\rm(Laufer, [\rLf])}
	 Let $M$ be a Stein \mfd, $D\subset M$ open, and
	$H=H^{p,q}(D)$, $p\ge0$, $q\ge1$, a Dolbeault cohomology group of $D$.
	 Then either $H=0$ or\/ $\dim_{\CC}H=\infty$.
}

	 In the remainder of this section we adopt the following.
	 Let $M$ a \cpx\ \bmfd, $D\subset M$ open, and $H=H^{p,q}(D)$, $p\ge0$, $q\ge1$,
	a Dolbeault cohomology group, or $H=H^q(D,\OO^{\La_p})$ a sheaf cohomology group
	with values in the sheaf $\OO^{\La_p}$ of germs of \holo\ sections of the \bvbdl\
	$\La_p\to M$ of $(p,0)$-forms.
	 (In \fdms\ the above-mentioned Dolbeault and sheaf cohomology groups are
	naturally isomorphic by the Dolbeault isomorphism \th.
	 In \idms\ the analog of the Dolbeault isomorphism is not yet proved except in
	very special cases, and sometimes may in fact fail.)
	 If $f\in\OO(D)$, $v\in\OO(D,T^{1,0}D)$, then $f$ and $v$ both act naturally
	on $H$ as linear operators.
	 We can set up these actions of multiplication and Lie derivation as follows.
	 We take the case of the Dolbeault group; the argument for the sheaf cohomology
	group is similar, only simpler.
	 If $\aa$ is a smooth $(p,q)$-form on $D$, then the product $f\aa$ and the
	Lie derivative $\LL_v\aa$  are smooth
	$(p,q)$-forms on $D$, and the commutation relations
	$\delbar(f\aa)=f\delbar\aa$ for the product and 
	$\delbar\LL_v\aa=\LL_v\delbar\aa$ for the Lie derivative
	show that the actions $M_f$ of multiplication by $f$ and $\LL_v$ of
	Lie derivation by $v$ descend to \cpx\ linear operators $[M_f],[\LL_v]\cln H\to H$.
	 As $\LL_v(f\aa)-f\LL_v\aa=(vf)\aa$, the commutation relation 
	$[\LL_v,M_f]=\LL_v M_f-M_f\LL_v=M_{vf}$ shows that
	$[[\LL_v],[M_f]]=\ad_{[\LL_v]}[M_f]=[M_{vf}]$.
	 Put $[f]=[M_f]$, $[v]=[\LL_v]$ for short.
	
	 Let $I\subset\OO(M)$ be the kernel of the representation $\OO(M)\to E=\End(H)$
	given by $f\mapsto[f|D]$; hence $f\in I$ if and only if $[f|D]=0$.
	 Henceforth we drop the restrictions from $M$ to $D$ from the notation.
	 Note that if $v\in\OO(M,T^{1,0}M)$ and $f\in I$, then 
	$[vf]=\ad_{[v]}[f]=\ad_{[v]}0=0$,
	i.e., $vf\in I$, so $I$ is closed under Lie derivation.
	 We call $I$ the \kiemel{annihilator (ideal)} of the group $H$ \wrt\ $M$.

\tetel{\t\tFB.}{Let $X$ be a \bspc, $\Oa\subset X$ \pscx\ open, $M\subset\Oa$
	a closed split \cpx\ \bsmfd\ of\/ $\Oa$, $D\subset M$ open, and $H$ a 
	cohomology group and $I$ its annihilator ideal as above.
	 If\/ {\rm(a)}, {\rm(b)}, or\/ {\rm(c)} below holds, then either $H=0$ or\/
        $\dim_\CC H=\infty$.
\vskip0pt
	{\rm(a)} $D$ (or $M$) admits a \recpr\ $f,v$.
\vskip0pt
	{\rm(b)} $X$ is \idml\ and has a \sbs, \pshdom\ is possible in every \pscx\ open 
	subset of\/ $\Oa$, and $M\subset\Oa$ is of a finite codimension $k\ge1$.
\vskip0pt
	{\rm(c)} $X$ has an \ubs, $M$ is connected, and $I$ has a common zero or $I$
	is sequentially closed.
}

\biz{Proof.}{Suppose for a contradiction that $1\le n=\dim_\CC H<\infty$.
	 If $I=0$, then the \idml\ \vspc\ $\OO(M)$ (or $\OO(D)$) is embedded in the \fdml\
	\vspc\ $E=\End(H)$ by the injective representation $\OO(M)\to E$
	induced by multiplication $f\mapsto[f|D]$.
	 Hence $I\not=0$.
	 Below in each case (a), (b), (c) we find a contradiction or 
	we show that $I=(1)$.
	 Then $H=1H=0$ is a contradiction that proves \t\tFB.

	 As mentioned in \p\tDF\ we may and do choose a basis of $H$ so that each matrix
	$[f]$ is upper triangular, and denote its diagonal characters 
	by $\chi_i{:}\;\OO(M)\to\CC$ for $i=1,\ldots,n$.
	 Let $K$ be the joint kernel of these characters, i.e., 
	$K=\{f\in\OO(M){:}\;\chi_i(f)=0, i=1,\ldots,n\}$.

	 The Cayley--Hamilton \th\ of linear algebra tells us that for every matrix 
	$[f]\in E$, $f\in\OO(M)$,
	there is a one-variable polynomial $p\in\CC[z]$ monic of degree $n$ \st\
	$0=p([f])=[p(f)]$, i.e., $p(f)\in I$.

	 Similarly, for every linear operator $\ad_{[v]}\in\End(E)$ there is a one-variable
	polynomial $p\in\CC[z]$ monic of degree $m=n^2$ \st\ $0=p(\ad_{[v]})$, i.e.,
	if $f\in\OO(M)$, then $0=p(\ad_{[v]})[f]=[p(v)f]$.
	 So $p(v)f=p(\LL_v)f$ belongs to the annihilator $I$ for all $f\in\OO(M)$, where
	$p(\LL_v)f$ is the resulting \fn\ obtained by applying the differential operator
	$p(\LL_v)$ to the \fn\ $f$.

	(a) Somewhat more generally, (a) is valid if $D$ is any \cpx\ \bmfd\ with
	a \recpr\ $f,v$.
	 Indeed, there is a one-variable polynomial $p\in\CC[z]$ monic of degree $n$ \st\
	$p(f)\in I$.
	 As $I$ is closed under Lie derivation, $1=\frac{1}{n!}v^np(f)\in I$ as well.

	(b) Consider the intersection $K\cap X^*$, where $X^*$ is the dual space of $X$.
	 The \vspc\ $K\cap X^*=\{\xi\in X^*{:}\;\chi_i(\xi_i|M)=0,i=1,\ldots,n\}$ is of
	codimension at most $n$ in the \idml\ \vspc\ $X^*$, hence $K\cap X^*$ itself is
	\idml.
	 There are $k+1$ linearly independent $\xi_\ka\in K\cap X^*$ for $\ka=0,\ldots,k$.
	 Letting $f_\ka=\xi_\ka|M$ we find by \p\tCG\ that $f_0,\ldots,f_k$ have no common
	critical points in $M$.
	 Moreover, $f_\ka^n$ belongs to $I$ for $\ka=0,\ldots,k$.
	 \t\tCX\ gives us \vfld{}s $v_\ka\in\OO(M,TM)$ for $\ka=0,\ldots,k$ \st\
	$\sum_{\ka=0}^k(v_\ka f_\ka)(x)=1$ for all $x\in M$.
	 Let $p_\ka\in\CC[z]$ be a one-variable polynomial monic of degree $m=n^2$ \st\
	$p(\ad_{[v_\ka]})=0$, and define $g_\ka\in\OO(M)$ by 
	$g_\ka=p_\ka(v_\ka)f_\ka^m=v_\ka^mf_\ka^m+\ldots$ for $\ka=0,\ldots,k$.
	 Note that $g_\ka$ belongs to $I$ for all $\ka$, and if $x_0\in M$ and 
	$f_\ka(x_0)=0$, then $g_\ka(x_0)=m!(v_\ka f_\ka)(x_0)^m$.

	 Consider the \fns\ $f_\ka^n,g_\ka$ for $\ka=0,\ldots,k$.
	 They all belong to $I$, and we claim that they have no common zeros in $M$.
	 Indeed, let $x_0\in M$ be any point and suppose that $f_\ka(x_0)^n=g_\ka(x_0)=0$
	for all $\ka=0,\ldots,k$.
	 Then $f_\ka(x_0)=0$, and so the equality $g_\ka(x_0)=m!(v_\ka f_\ka)(x_0)^m=0$
	implies that $(v_\ka f_\ka)(x_0)=0$.
	 Now $\sum_{\ka=0}^k(v_\ka f_\ka)(x_0)$ equals both $0$ and $1$; thus our
	members $f_\ka^n,g_\ka$, $\ka=0,\ldots,k$, of the ideal $I$ have no common zeros
	in $M$.
	 \t\tBB(d) gives us $a_\ka,b_\ka\in\OO(M)$ for $\ka=0,\ldots,k$
	with $1=\sum_{\ka=0}^k(a_\ka f_\ka^n+b_\ka g_\ka)\in I$ on $M$.

	(c) If $\chi_i{:}\;\OO(M)\to\CC$ is a \cts\ character, then
	$I$ has a common zero by \p\tDF.
	 (Nobody has ever seen a dis\cts\ character of an $\OO(M)$; 
	see Michael's problem in [\rMc] or [\rM].)
	 Suppose that $x_0\in M$ is a common zero of $I$.
	 As $I\not=0$ and $M$ is connected, there is an $f\in I$ with vanishing order
	$1\le s<\infty$ at $x_0$.
	 Nonzero is the Fr\'echet derivative $f^{(s)}(x_0)$ relative to a biholomorphism 
	of an open \nbd\ $U\subset M$ of $x_0$ onto an open subset $V$ of a \bspc.
	 There are vectors $\xi_1,\ldots,\xi_s\in T_{x_0}M$ with 
	$f^{(s)}(x_0)\xi_1\ldots\xi_s=1$.
	 There are \vfld{}s $v_i\in\OO(M,TM)$ with $v_i(x_0)=\xi_i$ for $i=1,\ldots,s$, e.g.,
	of the form $v_i(x)=P'(x)\xi_i$, where $P'\in\OO(M,\End(X))$ is as in the proof of
	\t\tED.
	 The \fn\ $g=\LL_{v_1}\ldots\LL_{v_s}f\in\OO(M)$ belongs to $I$, hence $g(x_0)=0$
	on the one hand.
	 On the other hand $g(x_0)=f^{(s)}(x_0)\xi_1\ldots\xi_s=1$ by \p\tEA; a contradiction.

	 If $I$ is sequentially closed and has no common zeros in $M$, then $1\in I$ 
	by \t\tED.
	 The proof of \t\tFB\ is complete.
	 
}

	 \t\tFB\ applies to all finite codimensional closed \cpx\ \hsmfd{}s of $\ell_2$,
	and also to some \idml\ ones such as those in \t\tCC(b) (mapping spaces).
	 It seems likely (but currently unknown) that the conclusion of \t\tFB\ also
	holds for every closed \cpx\ \hsmfd\ $M$ of $\ell_2$, because $M$ might
	always admit a \recpr\ $f,v$.
	 It is already known, see [\rPt], that there are a nowhere zero \vfld\
	$v\in\OO(M,TM)$ and a nowhere zero covector field $\oa\in\OO(M,T^*M)$, but
	it seems unknown whether such an $\oa$ can be chosen to be exact or even closed
	on $M$.

	 The proof of \t\tFB\ has much in common with that of Laufer's for \t\tFA,
	but unlike in \fdms\ finitely many numerical \fns\ do not separate the points of an
	\idml\ \bmfd.
	 To overcome this difficulty, we work harder with the vector fields $v$.

\alcim{\sG. NOWHERE CRITICAL HOLOMORPHIC FUNCTIONS.}

	 In this section we look at a simple mechanism that extends 
	the special case $X=\ell_1$ of \t\tCD\ to all
	closed \idml\ split \cpx\ \bsmfd{}s of $\ell_1$.

\tetel{\p\tGA.}{Let $T$ be a separable \tspc, $X,Y$ \bspc{}s, $Y$ nonseparable, $Z=X\times Y$
	with the product norm\/ $\|z\|=\max\{\|x\|,\|y\|\}$ for $z=(x,y)\in Z$, $H=\Hom(X,Y)$
	with the operator norm, and $A\in C(T,H)$ a \bdd\ \cts\ \fn\ with operator values.
	 Denote for a map $f{:}\;X\to Y$ its graph by\/ $\Ga(f)=\{(x,y)\in Z{:}\;y=f(x)\}$.
	 We then have for the set $E=\bigcup_{t\in T}\Ga(A(t))$ that
\vskip0pt
	{\rm(a)} $\overline{E}\cap(\{x_0\}\times Y)=F_{x_0}$ is separable for every $x_0\in X$, and
\vskip0pt
	{\rm(b)} the closure $\overline{E}$ is nowhere dense in $Z$.
}

\biz{Proof.}{As it implies (b) let us prove (a).
	 Let $T'\subset T$ be countable and dense in $T$, and define the countable subset
	$F'=\{(x_0,A(t)x_0){:}\;t\in T'\}$ of $F=F_{x_0}$.
	 As $z_0=(x_0,y_0)\in\overline{E}$ there is a sequence $z_n=(x_n,A(t_n)x_n)\in\Ga(A(t_n))$,
	$t_n\in T$, $n\ge1$, with $z_n\to z_0$ as $n\to\infty$; in particular, $\lim_{n\to\infty}x_n=x_0$.
	 Our $A$ being \cts\ at each $t_n$, $n\ge1$, there is a $t'_n\in T'$ so close to $t_n$ in $T$
	that $\|A(t_n)-A(t'_n)\|<1/n$.
	 The bound $M=\sup_{t\in T}\|A(t)\|$ is finite by assumption.
	 Look at $z'_n=(x_0,A(t'_n)x_0)\in F'$ for $n\ge1$.
	 As $\|z_n-z_0\|=\max\{\|x_n-x_0\|,\|A(t_n)x_n-A(t'_n)x_0\|\}\le
	\|x_n-x_0\|+\|A(t_n)x_n-A(t_n)x_0\|+\|A(t_n)x_0-A(t'_n)x_0\|\le
	\|x_n-x_0\|+M\|x_n-x_0\|+\|x_0\|/n\to0$ as $n\to\infty$, we find that $\lim_{n\to\infty}z'_n=z_0$,
	i.e., the countable subset $F'$ of $F$ is dense in $F$.
	 The proof of \p\tGA\ is complete.
}

\tetel{\p\tGB.}{Let $X,Y$ be \bspc{}s, $Z=X\times Y$, $\Oa\subset X$ open, $m{:}\;\Oa\to Y$
	with \cts\ Fr\'echet derivative $m'$ \bdd\ on\/ $\Oa$, 
	and $M=\Ga(m)\subset\Oa\times Y$ the graph of $m$.
	 If $X$ is separable but its dual space $X^*$ is nonseparable, then there is a
	nowhere dense closed subset $B\subset Z^*$ \st\ for $f\in Z^*\setminus B$ the restriction $f|M$
	is nowhere critical on $M$.
}

\biz{Proof.}{As $Z^*=X^*\times Y^*$, we may write any $f\in Z^*$ uniquely as
	$f(x,y)=\xi(x)+\eta(y)$ with $\xi\in X^*$ and $\eta\in Y^*$.
	 A point $z_0=(x_0,y_0)$ with $y_0=m(x_0)$ is a critical point of $f|M$ if and only if
	$x=x_0$ is a critical point of $f(x,m(x))$, i.e., $\xi+\eta m'(x_0)=0$ in $Z^*$.
	 In other words, $f|M$ is critical at $z_0$ if and only if $\xi=A(x_0)\eta$, where
	$A(x)\in\Hom(Y^*,X^*)$ is given for $x\in\Oa$ 
	by the transpose of $-m'(x)\in\Hom(X,Y)$.
	 The set of `exceptional' \fns\ $f$ with $f|M$ critical at some point of $M$ is thus
	of the form $E=\bigcup_{x\in\Oa}\Ga(A(x))$ in $Z^*=X^*\times Y^*$.
	 \p\tGA\ applies and shows that $B=\overline{E}$ is nowhere dense closed and does the job.
}

\tetel{\p\tGC.}{Let $Z$ be a separable \bspc, $\Oa\subset Z$ open, and $M$ a closed $C^1$-smooth
	split \bsmfd\ of\/ $\Oa$.
	 If the cotangent space $T^*_xM$ is nonseparable for all $x\in M$, then there is a
	dense $G_\da$ subset $H$ of the dual space $Z^*$ of $Z$ \st\ of each $f\in H$ the
	restriction $f|M$ is nowhere critical on $M$.
}

\biz{Proof.}{As locally the Lindel\"of space $M$ is given by a graph $y=m(x)$
	of a \fn\ $m$ \bdd\ in the $C^1$-norm, our claim follows from \p\tGB\ 
	and the Baire category \th\ applied to the \bspc\ $Z^*$.
}

\tetel{\t\tGD.}{Let $\Oa\subset\ell_1$ be open and $M$ a closed split \cpx\ \bsmfd\ of\/ $\Oa$.
	 If $M$ is \idml\ at each of its points, then $M$ admits a nowhere critical \holo\ \fn\
	$f\in\OO(M)$ that can be extended to a nowhere critical
	\holo\ \fn\ $\tilde f\in\OO(\ell_1)$ on $\ell_1$.
	 In fact, $\tilde f$ can be taken linear and to be an arbitrary member of 
	a dense $G_\da$ subset of the dual space $\ell_1^*=\ell_\infty$.
	 Further, if\/ $\Oa$ is \pscx, then there is a $v\in\OO(M,T^{1,0}M)$ 
	with $vf=1$ on $M$.
}

\biz{Proof.}{It is a famous \th\ of Pe{\l}czy{\'n}ski's, see [\rP] or [\rLT, Thm.\,2.a.3],
	that any closed complemented \idml\ linear subspace 
	(such as the tangent space $T_xM$ for $x\in M$) of $\ell_1$ is isomorphic to $\ell_1$,
	hence it is separable but its dual (such as the cotangent space $T^*_xM$ for $x\in M$)
	is not, being isomorphic to $\ell_\infty$.
	 \p\tGC\ yields the $f$, and \t\tCX\ a $v$ for $f$, completing the proof of \t\tGD.
}

         It would be interesting to know whether one could use in a similar manner
        \bspc{}s of polynomials of higher degree.
         Note that many of the standard separable \bspc{}s $Z$ (e.g., $Z=\ell_2$)
        admit nonseparable \bspc{}s of polynomials over them (e.g., quadratic forms
        over $\ell_2$) even if their duals $Z^*$ themselves are still separable.
         Another question is whether a \cpx\ algebraic submanifold $M$ of $\CC^n$ admits
        a \holo\ polynomial $f$ nowhere critical on $\CC^n$ whose restriction $f|M$ is
        also nowhere critical on $M$.

        {\vekony Acknowledgement.} The essential idea of this paper occured to the
        author at one of the annual summer workshops on analysis and probability at Texas
        A\,\&\,M University.
         He is grateful for their hospitality, and to L\'aszl\'o Lempert for a remark
        that led to the above \S\,7.

\vskip0.05truein
\centerline{\scVIII References}
\vskip0.05truein
\baselineskip=11pt
\parskip=0pt
\frenchspacing
{\rmVIII

	[\rD] Dineen, S.,
	{\itVIII 
	Complex analysis in infinite dimensional spaces},
	Springer-Verlag, 
	London,
	(1999).

	[\rDPV] \vonal, Patyi, I., Venkova, M.,
	{\itVIII 
	Inverses depending holomorphically on a parameter in a Banach space},
	J. Funct. Anal., 
	{\bfVIII 237} (2006), no. 1, 338--349.

	[\rF] Forstneri\v{c}, F.,
	{\itVIII
	Noncritical holomorphic functions on Stein manifolds},
	Acta Math.,
	{\bfVIII 191} (2003), 143--189.

	[\rLf] Laufer, H.\,B.,
	{\itVIII
	On the infinite dimensionality of the Dolbeault cohomology groups},
	Proc. Amer. Math. Soc.,
	{\bfVIII 52} No. 1 (1975), 293--296.

	[\rLA] Lempert,~L.,
	{\itVIII
	The Dolbeault complex in infinite dimensions~I},
	J. Amer. Math. Soc., {\bfVIII 11} (1998), 485--520.

	[\rLB] \vonal,
	{\itVIII
	Plurisubharmonic domination},
	J. Amer. Math. Soc.,
	{\bfVIII 17}
	(2004),
	361--372.

	[\rLP] \vonal, Patyi,~I.,
	{\itVIII
	Analytic sheaves in Banach spaces},
	Ann. Sci. \'Ecole Norm. Sup., 
	S\'er. 4,
	{\bfVIII 40}
	(2007),
	453--486.

	[\rLT] Lindenstrauss, J., Tzafriri, L.,
	{\itVIII
	Classical Banach Spaces I, Sequence Spaces},
	Spring\-er, Berlin, (1977).

	[\rMc] Michael, E.\,A.
	{\itVIII
	Locally multiplicatively-convex topological algebras},
	Mem. Amer. Math. Soc., {\bfVIII 1952}, (1952). no. 11, 79 pp. 

	[\rM] Mujica, J.,
	{\itVIII
	 Complex analysis in Banach spaces},
	North--Holland, Amsterdam,
	(1986).

        [\rPt] Patyi, I.,
        {\itVIII
        On holomorphic Banach vector bundles over Banach spaces},
	Math. Ann.,
	{\bfVIII 341} (2008), no. 2, 455--482.

	[\rP] {Pe{\l}czy{\'n}ski, A.},
     	{\itVIII
	Projections in certain {B}anach spaces},
   	Studia Math.,
    	{\bfVIII 19}
      	(1960),
    	209--228.

	[\rS] Schottenloher, M.,
	{\itVIII
	Spectrum and envelope of holomorphy for infinite-dimensional Riemann domains},
	Math. Ann.,
	{\bfVIII 263} (1983), no. 2, 213--219. 

	[\rZ] Zerhusen,~A.\,B.,
	{\itVIII
	Embeddings of pseudoconvex domains in certain Banach spaces},
	Math. Ann.,
	{\bfVIII 336}
	(2006), no. 2,
	269--280.

}
\vskip0.10truein
\centerline{\vastag*~***~*}
\vskip0.10truein
{\scVIII
        Imre Patyi,
        Department of Mathematics and Statistics,
        Georgia State University,
        Atlanta, GA 30303-3083, USA,
        {\ttVIII matixp\@langate.gsu.edu}
}
\bye